\documentclass{amsart}

\usepackage{amssymb}
\usepackage{graphicx} 
\usepackage{eepic}
\usepackage{amsmath}
\usepackage{caption}
\usepackage{graphs}
\usepackage{subfigure}
\usepackage{color}

\definecolor{red}{rgb}{0.9,0.1,0.1}
\newtheorem{theorem}{Theorem}[section]
\newtheorem{lemma}[theorem]{Lemma}
\newtheorem{proposition}[theorem]{Proposition}

\theoremstyle{definition}

\newtheorem{corollary}[theorem]{Corollary}
\newtheorem{example}[theorem]{Example}

\newtheorem{conjecture}[theorem]{Conjecture}

\theoremstyle{remark}

\numberwithin{equation}{section}

\newcommand{\rank}[1]{\operatorname*{rank\left(#1\right)}}
\newcommand{\M}{\mathcal{M}}
\newcommand{\X}{\mathbf{X}}

\newcommand{\sh}[1]{\medskip \noindent {\bf #1}}
\newcommand{\field}{{\bf k}} 
\newcommand{\Field}{{\bf K}} 

\DeclareMathOperator{\In}{In}
\DeclareMathOperator{\lk}{link}
\newcommand{\PG}{\mathcal{PG}}

\begin{document}
\bibliographystyle{abbrv}
\title[]{Face Numbers of Certain Cohen-Macaulay Flag Complexes}
\author{Jonathan Browder}
\date{10 October 2010}

\begin{abstract}
We show that if a $d$-dimensional Cohen-Macaulay complex is, in a certain sense, sufficiently ``close'' to being balanced, then there is a $d$-dimensional balanced Cohen-Macaulay complex having the same $f$-vector. This in turn provides some partial evidence for a conjecture of Kalai on the $f$-vectors of Cohen-Macaulay flag complexes.

\end{abstract}

\maketitle
 
\section{Introduction} 
\subsection{Background} One of the fundamental invariants of a simplicial complex $\Delta$ is its $f$-vector, $f(\Delta) = (f_{-1}, f_0, \ldots , f_{\dim(\Delta)})$, which lists the number of faces $\Delta$ has in each dimension (i.e., $f_i$ is the number of $i$-dimensional faces of $\Delta$). Characterizing the possible $f$-vectors of various classes of simplicial complexes is one of the central problems of geometric combinatorics. Of particular recent interest are flag complexes and balanced complexes; it was conjectured by Kalai and proven by Frohmader \cite{MR2391144} that the $f$-vector of an arbitrary flag complex is also the $f$-vector of some balanced complex (though the reverse does not hold). Kalai further made the following conjecture, which remains open:

\begin{conjecture} \label{kalaiconj} Let $\Delta$ be a \emph{Cohen-Macaulay} flag complex of dimension $d-1$. Then there is a $(d-1)$-dimensional Cohen-Macaulay balanced complex $\Gamma$ such that $f(\Delta) = f(\Gamma)$.
\end{conjecture}

Our main theorem provides some partial evidence for this conjecture.

\begin{theorem}\label{main} Let $\Gamma_1, \Gamma_2, \ldots, \Gamma_k$ be $0$ or $1$ dimensional flag complexes (i.e., triangle-free graphs) such that for each $i$, either $\Gamma_i$ is bipartite or  $\Gamma_i - e$ is bipartite for some edge $e$ of $\Gamma_i$. Let $\Gamma = \Gamma_1 * \Gamma_2 * \cdots * \Gamma_k$ (where $*$ denotes the simplicial join). Then for $\Delta$ any full-dimensional Cohen-Macaulay subcomplex of $\Gamma$, $f(\Delta)$ is the $f$-vector of some balanced complex of the same dimension.
\end{theorem}

Notice that the complexes described in Theorem \ref{main} are in some sense ``close'' to balanced; they can be made balanced by deleting an appropriate edge from each of the terms in the join which are not bipartite. In Section \ref{examples} we will see that the theorem applies to a large class of examples of flag complexes arising as independence complexes of graphs with certain properties. Note, however, that the theorem applies to many complexes which are not flag; while the complex $\Gamma$ is flag, the subcomplexes described need not be.

\subsection{Preliminaries}We begin by reviewing some basic concepts and notation from the study of simplicial complexes. For further details, \cite{MR1453579} is a good reference.

\sh{Simplicial Complexes and Multicomplexes} A \emph{simplicial complex} $\Delta$ on finite vertex set $V$ is a set of subsets of $V$ which is closed under inclusion. An element of $\Delta$ is called a \emph{face}, the faces which are maximal with respect to inclusion are called \emph{facets}.  The dimension of a face $\gamma$ is $\dim(\gamma) := |\gamma| - 1$, and the dimension of the complex is $\dim(\Delta) := \max\{ \dim(\tau) : \tau \in \Delta \}$. Faces of dimension $0$ and $1$ are called vertices and edges, respectively. The complex is \emph{pure} if all of its facets have the same dimension. 
For $i \leq \dim(\Delta)$, the \emph{$i$-skeleton} of $\Delta$ is the subcomplex of $\Delta$ consisting of all the faces of $\Delta$ with dimension no greater than $i$. In particular, the 1-skeleton of $\Delta$ may be thought of as a graph. For $\tau \in \Delta$, the \emph{link} of $\tau$ in $\Delta$ is $\lk_{\Delta}(\tau) := \{ \gamma \in \Delta : \tau \cap \gamma = \emptyset, \gamma \cup \tau \in \Delta \}$.

The $f$-vector of simplicial complex $\Delta$ is defined to be $f(\Delta) = (f_{-1}, f_0, \ldots , f_{d-1})$, where $d-1$ is the dimension of $\Delta$ and $f_i$ is the number of $i$-dimensional faces of $\Delta$ (these $f_i$ are known as the \emph{face numbers} of $\Delta$). Notice that $f_{-1} = 1$ for any non-empty $\Delta$, as the empty set will be the unique $(-1)$-dimensional face. 

In practice it is often more convenient to study the face numbers of the complex in terms of the $h$-vector of the complex, $h(\Delta) := (h_0, h_1, \ldots , h_d)$, where the numbers $h_i$ are defined by the relation
\begin{equation*}
\sum_{i = 0}^d h_ix^i = \sum_{i=0}^d f_{i-1}x^i(1-x)^{d - i}.
\end{equation*}
It is clear that the $f$-vector of $\Delta$ completely determines its $h$-vector and vice versa.

Similary, let $\X$ be a finite set of variables, and define a \emph{multicomplex} on $\X$ to be a collection of monomials in $\X$ which is closed under divisibility (we include $1$ as the unique degree $0$ element of any non-empty multicomplex). Notice that if $M$ is a  multicomplex on $\X$ such that every element of M is square-free, then $M$ corresponds to a simplicial complex in the obvious way. The $F$-vector of a multicomplex $M$ is $F(M) := (F_0, F_1, \ldots )$, where $F_i$ is the number of elements in $M$ of degree $i$ (if $M$ is also a simplicial complex then the $F$-vector is just the $f$-vector up to a shift in index). 

For $S \subseteq \X$ and $m$ a monomial on $\X$, let $m_S$ denote the part of $m$ supported in $S$ (i.e., the unique monomial such that $m = m_Sm_{X-S}$, where $m_{X-S}$ is divisible by no element of $S$).

\sh{Stanley-Reisner Rings and the Cohen-Macaulay Property.} Let $\Delta$ be a $(d-1)$-dimensional simplicial complex on vertex set $V$, and let $\X = \{ x_v : v \in V \}$ be a set of variables indexed by $V$. Fix a field $\field$ of characteristic zero, and let $\field[\X]$ be the polynomial ring over $\field$ in the variables of $\X$ (with the grading $\deg{x_v} = 1$). Then the \emph{Stanley-Reisner ring} of $\Delta$ is $\field[\Delta] := \frac{\field[\X]}{I_{\Delta}}$, where $I_{\Delta}$ is the ideal in $\field[\X]$ generated by the squarefree monomials $x_{v_1}x_{v_2}\ldots x_{v_k}$ such that $\{ v_1, v_2, \ldots v_k \} \notin \Delta$. We call $I_{\Delta}$ the Stanley-Reisner ideal of $\Delta$; it is easy to see that it is generated by the monomials corresponding to the minimal non-faces of $\Delta$.

The Krull dimension of $\field[\Delta]$ is $d$ \cite{MR1453579}. A \emph{linear system of parameters} (l.s.o.p.) for $\field[\Delta]$ is a sequence $\theta_1, \theta_2, \ldots, \theta_d$ of elements of $\field[\Delta]_1$ such that $\frac{\field[\Delta]}{(\theta_1, \theta_2, \ldots, \theta_d)}$ is finite-dimensional as a $\field$-vector space. It follows from Noether Normalization that some l.s.o.p. must exist.

We will define $\Delta$ to be $\field$-Cohen-Macaulay ($\field$-CM) if for some (equivalently, every) l.s.o.p. $\theta_1, \theta_2, \ldots, \theta_d$ of $\field[\Delta]$, 
$$
\dim_{\field} (\field[\Delta]/(\theta_1, \theta_2, \ldots, \theta_d))_i = h_i  \quad \mbox{for all } 0\leq i \leq d.
$$

When the field $\field$ is understood, we will simply say such a complex is Cohen-Macaulay (CM). Note that if $\Delta$ is $\field$-CM, then $\Delta$ is $\Field$-CM for each field $\Field$ with the same characteristic as $\field$.

There are many equivalent definitions of the Cohen-Macaulay property, see, for example \cite{MR1453579}. Of particular note is Reisner's characterization \cite{MR0407036} of CM complexes in terms of the vanishing of certain homologies. In particular, it follows from Reisner's result that every CM complex is pure. Many interesting complexes are CM, including all shellable complexes and all triangulations of balls and spheres.

\sh{Balanced Complexes and Flag Complexes} 
A simplicial complex $\Delta$ is \emph{flag} if every clique in the 1-skeleton of $\Delta$ forms a face of $\Delta$. In particular $\Delta$ is completely determined by its set of edges, and $I_{\Delta}$ is generated in degree two. In this case $\Delta$ is both the clique complex of its $1$-skeleton and the independence complex of the graph complement of its 1-skeleton.

For $\Delta$ a simplicial complex on $V$, a map $\kappa : V \rightarrow [k]$ is called a \emph{proper $k$-coloring} of $\Delta$ if whenever distinct vertices $v_1$ and  $v_2$ are contained in a common face of $\Delta$, $\kappa(v_1) \neq \kappa(v_2)$ (in other words, $\kappa$ is a proper coloring of the 1-skeleton of $\Delta$ in the graph theoretic sense). A complex which has a proper $k$-coloring is called \emph{$k$-colorable}.  If $\Delta$ is $(d-1)$-dimensional it is clear that a minimum of $d$ colors are needed for a proper coloring of $\Delta$; if $\Delta$ is in fact $d$-colorable we say $\Delta$ is \emph{balanced}. (Sometimes these complexes are called \emph{completely balanced} if more general types of balance are in play.) A result of Stanley \cite{MR563925} (necessity) and Bj\"orner, Frankl and Stanley \cite{MR905148} (sufficiency) completely characterized the $h$-vectors of balanced Cohen-Macaulay complexes.

\begin{theorem} \cite{MR905148} \label{ch} Let $h = (h_0, h_1, \ldots , h_d)$. The following are equivalent:
\begin{enumerate}
\item There exists a balanced $(d-1)$-dimensional Cohen-Macaulay complex $\Delta$ such that $h(\Delta) = h$.
\item There exists d-colorable simplicial complex $\Gamma$ such that $f(\Gamma) = h$.
\end{enumerate}

\end{theorem}

It is furthermore worth noting that a purely numerical characterization of the $f$-vectors of d-colorable simplicial complexes was found in \cite{MR1018807}; if Conjecture \ref{kalaiconj} is true, it would imply that the $h$-vector of any Cohen-Macaulay flag complex is the $f$-vector of such a complex.

\begin{example} Consider the simplicial complex $\Gamma$ shown in Figure \ref{fig1}. Clearly $\Gamma$ is $3$-colorable, so there is a balanced 2-dimensional Cohen-Macaulay complex $\Delta$ such that $h(\Delta) = f(\Gamma) = (1,4,5,1)$. Then $f(\Delta) = (1,7,16,11)$.

Now suppose there is some 2-dimensional flag complex $\Omega$ with $f(\Omega) = f(\Delta)$. The $1$-skeleton of $\Omega$ is then a graph on $7$ vertices which contains no $K_4$ (as $\Omega$ is 2-dimensional), and has 16 edges. Tur\'an's Theorem \cite{MR0018405} tells us this is the maximum number of edges in a $K_4$-free graph on $7$ vertices, and in particular that the $1$-skeleton of $\Omega$ must in fact be the Tur\'an graph $T(7,3)$. But $T(7,3)$ contains 12 triangles, all of which must be faces of $\Omega$, a contradiction. Thus there is no 2-dimensional flag complex having the same $f$-vector as $\Delta$, so the reverse of Conjecture \ref{kalaiconj} does not hold.

\begin{figure}[h] 
 \includegraphics[scale=0.2]{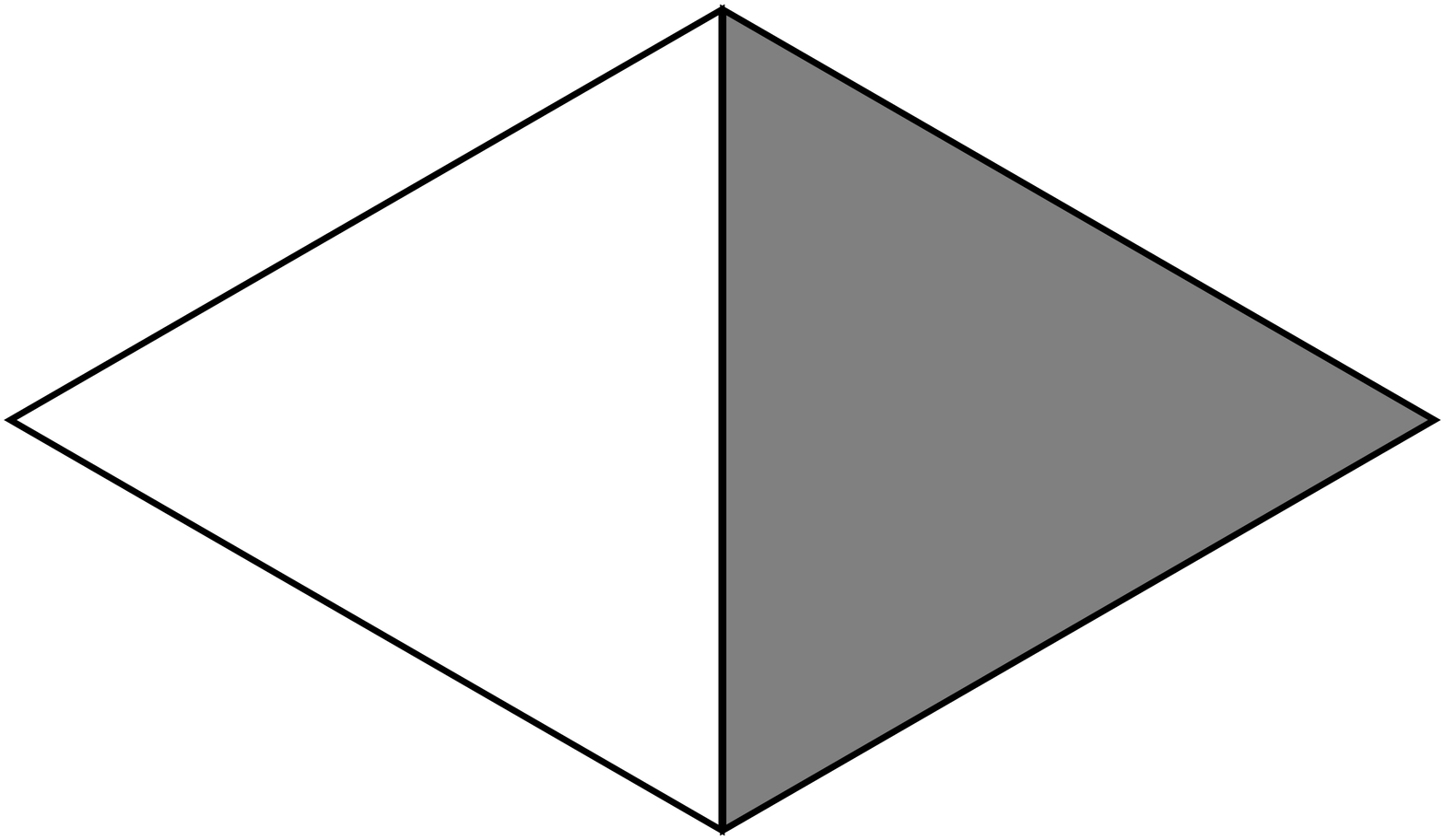}
 \caption{}
 \label{fig1}
\end{figure}
\end{example}

\section{Proof of the Main Theorem}
\subsection{Outline}We first outline our approach, which is adapted from that used in \cite{MR2122283} and \cite{genbal}. Throughout the following, let $\Delta$ be a $(d-1)$-dimensional $\field$-CM complex on $V$, $n = |V| $,  and $\X = \{ x_v : v \in V \}$.

Suppose we fix a total order on $\X$ and let $\prec$ denote the corresponding reverse lexicographical (revlex) order on the monomials in $\X$. Let $T_{\prec}$ denote the last $d$ elements of $X$ with respect to $\prec$, and $(T_{\prec})$ the ideal in $\field[\X]$ these elements generate. Further suppose we pick a graded automorphism $g$ of $\field[\X]$ (considered as a matrix with entries in $\field$) such that $T_{\prec}$ is an l.s.o.p. of $\frac{\field[\X]}{gI_{\Delta}}$. Then $\dim_k \left(\frac{\field[\X]}{gI_{\Delta}+ (T_{\prec})} \right)_i = h_i$, as $\frac{\field[\X]}{gI_{\Delta}}$ is isomorphic to $\field[\Delta]$ and $\Delta$ is CM. Furthermore, $gI_{\Delta}+ (T_{\prec})$ is a homogeneous ideal, so its revlex initial ideal $\In(gI_{\Delta}+ (T_{\prec}))$ is well-defined, and Theorem 15.3 of \cite{MR1322960} asserts that the set $B_g(\Delta)$ of monomials in $\X$ \emph{not} in $\In(gI_{\Delta}+ (T_{\prec}))$ is a $\field$-basis for $\frac{\field[\X]}{gI_{\Delta}+ (T_{\prec})}$. Then $B_g(\Delta)$ is a multicomplex (as $\In(gI_{\Delta}+ (T_{\prec}))$ is an ideal) with $F_i(B_g(\Delta)) = h_i$. Note that $T_{\prec} \subseteq \In(gI_{\Delta}+ (T_{\prec}))$, so $B_g(\Delta)$ is a multicomplex on $\X - T_{\prec}$.

Thus, in light of Theorem \ref{ch}, to prove Theorem \ref{main} it will suffice to show that for $\Delta$ a complex as in the statement of the theorem, we may choose an order on $\X$, automorphism $g$, and partition of $X - T_{\prec}$ into disjoint sets  $X_1, X_2, \ldots,  X_d$ such that $T_{\prec}$ is an l.s.o.p for $\frac{\field[\X]}{gI_{\Delta}}$, and for each $m \in B_g(\Delta)$ and $1 \leq i \leq d$, $\deg(m_{X_i}) \leq 1$ (the last condition ensures that $B_g(\Delta)$ is a simplicial complex with $d$-coloring corresponding to the partition of $X-T_{\prec}$).

Finally, to verify that $T_{\prec}$ is an l.s.o.p for $\frac{\field[\X]}{gI_{\Delta}}$, it will suffice to check that $g$ satisfies the \emph{Kind-Kleinschmidt condition} \cite{MR534824}:
\begin{itemize}
\item For every facet $\{ x_{v_1}, x_{v_2}, \ldots , x_{v_k} \}$ of $\Delta$, the submatrix of $g^{-1}$ given by the intersection of the last $d$ columns of $g$ with the rows corresponding to $v_1, v_2, \ldots, v_k$ has rank $k$.
\end{itemize}

It will be convenient in our arguments to replace the field $\field$ with a larger field $\Field$ defined to be the field of rational functions over $\field$ in indeterminates $z_1, z_2, z_3, z_4$. A complex which is $\field$-CM is also $\Field$-CM, and passing between the two will not affect the enumerative consequences of our arguments.

\subsection{Proof of Theorem \ref{main}} Let $\Gamma$ be a pure $d-1$-dimensional complex on $V$, $\X = \{ x_v : v \in V \}$. Let $\prec$ be a total order of $\X$, $g$ a graded automorphism of $\Field[\X]$. We will call $(\prec, g)$ a \emph{balancing pair} for $\Delta$ if there exists a partition of $X - T$ into disjoint sets  $X_1, X_2, \ldots,  X_d$ such that
\begin{enumerate}
\item $g$ satisfies the Kind-Kleinschmidt condtion for $\Delta$.
\item If $m \in B_g(\Delta)$ and $1 \leq i \leq d$, $\deg(m_{X_i}) \leq 1$.
\end{enumerate}

Then to prove Theorem \ref{main}, it suffices to show that there exist balancing pairs for the complexes in question. We will induct from the $k=1$ case with the aid of some lemmas.

First suppose we have balancing pairs $(g_1, \prec_1)$ and $(g_2, \prec_2)$ for complexes $\Delta_1$ and $\Delta_2$, respectively, where $\Delta_i$ has dimension $d_i-1$, vertex set $V_i$, and corresponding set of variables $\X^i$. Let $\X^i_1, \ldots \X^i_{d_i}$ denote the corresponding partitions of $\X^i - T_{\prec_i}$, for $i=1,2$.

Further suppose that for $i=1,2$,  $g_i$ is of the form
$$g_i = 
               \left[ \begin{array}{cc}
                A_i & \begin{array}{c} B_i\end{array}\\
                0 & C_i \end{array} \right], $$
where $A_i$ is a $(|V_i| - d_i) \times (|V_i|-d_i)$ matrix, $B_i$ a $(|V_i|-d_i) \times d_i$ matrix, and $C_i$ a $d_i \times d_i$ matrix. 

Now let $\prec$ be the order of $\X^1 \cup \X^2$ given by $x \prec y$ if and only if either
\begin{itemize}
\item$x, y \in \X^i$, and $x \prec_i y$,
\item $x \notin T_{\prec_i}$ for $i = 1, 2$ and $y \in T_{\prec_i}$ for some $i$, 
\item $x \in T_{\prec_1}$, $y \in T_{\prec_2}$, or
\item $x \in \X^1 - T_{\prec_1}$, $y \in \X^2 - T_{\prec_2}$.
\end{itemize}
Finally let
$$g =
          \left[ \begin{array}{cccc}
                 A_1 & 0 & B_1 & 0\\
                 0 & A_2 & 0 & B_2\\
                 0 & 0& C_1&0\\
                 0& 0& 0& C_2\end{array}\right].$$

\begin{lemma}\label{induct} The pair $(\prec, g)$ defined above is a balancing pair for $\Delta_1 * \Delta_2$.
\begin{proof} 
Let $\X = \X_1 \cup \X_2$. It is clear that $g$ is a graded automorphism of $\Field[\X]$, and observe that for $i=1,2,$
$$g_i^{-1} = 
               \left[ \begin{array}{cc}
                A_i^{-1} & \begin{array}{c} D_i\end{array}\\
                0 & C_i^{-1} \end{array} \right], $$
for some $D_i$, so
$$g^{-1} =
          \left[ \begin{array}{cccc}
                 A_1^{-1} & 0 & D_1 & 0\\
                 0 & A_2^{-1} & 0 & D_2\\
                 0 & 0& C_1^{-1}&0\\
                 0& 0& 0& C_2^{-1}\end{array}\right].$$

The dimension of $\Delta_1 * \Delta_2$ is $d_1 + d_2 -1$, and any facet $\tau$ of $\Delta_1 * \Delta_2$ is of the form $\tau = \tau_1 \cup \tau_2$, where $\tau_i$ is a facet of $\Delta_i$. Then the submatrix of $g^-1$ given by the intersection of the last $d_1 + d_2$ columns of $g$ with the rows indexed by $\tau$ is just
$$ M = \left[ \begin{array}{cc}
               M_{1,1} &0 \\ 0 &M_{2,1}\\
               M_{1,2} &0\\
               0 &M_{2,2} 
               \end{array} \right],
 $$
where for $i=1,2$, 
$$M_i = \left[ \begin{array}{c} M_{i,1}\\M_{i,2}\end{array}\right]$$ 
is the submatrix of $g_i^{-1}$ given by the intersection of the last $d_i$ columns of $g_i^{-1}$ with the rows indexed by $\tau_i$. Thus the rank of $M$ is $\rank{M_1} + \rank{M_2} = |\tau_1| + |\tau_2| = |\tau|$, so $g$ satisfies the Kind-Kleinschmidt condition for $\Delta_1 * \Delta_2$.

For our partition of $\X - T_{\prec}$, we will simply use that inherited from the partitions of $X_1 - T_{\prec_1}$ and $X_2 - T_{\prec_2}$ (noticing that $T_{\prec} = T_{\prec_1} \cup T_{\prec_2})$.

Suppose that for some $j=1,2$ and $i \in \{1, 2, \ldots, d_j \}$, $m$ is a monomial on $\X^j_i$ of degree greater than $1$. Then $m \notin B_{g_j}(\Delta_j)$, so there is some $\nu \in I_{\Delta_j}$ such that $m = \In(g_j\nu)$. But then it is clear that $\nu \in I_{\Delta}$, and $g\nu = g_j\nu$. As $g_j\nu$ involves only variables in $\X^j$ and $\prec$ restricts to $\prec_j$ on $\X^j$, $\In(g\nu) = m$, so $m\notin B_{\prec}$. Then as $B$ is a multicomplex, no monomial on in $\Field[\X]$ has degree greater than $1$ in $\X^j_i$.

\end{proof}
\end{lemma}

\begin{lemma} Suppose $\Delta$ is a full-dimensional subcomplex of $\Gamma$ and $(\prec, g)$ is a balancing pair for $\Gamma$. Then $(\prec, g)$ is a balancing pair for $\Delta$.

\begin{proof}As each face of $\Delta$ is a face of $\Gamma$, it following immediately that $g$ satisfies Kind-Kleinschmidt for $\Delta$.

Now, suppose $m$ is a monomial in $\Field[\X]$ such that for some $i \in \{1, 2, \ldots d \}$, $\deg(m_{X_i}) > 1$. Then $m \notin B_g(\Gamma)$, so $m \in \In(gI_{\Gamma})$. In other words, there is an element $\nu$ of $I_{\Gamma}$ such that $In(g\nu) = m$. But as $\Delta$ is a subcomplex of $\Gamma$, $I_{\Gamma} \subseteq I_{\Delta}$, so $m \in \In(gI_{\Delta})$.
\end{proof}
\end{lemma}

Theorem \ref{main} now follows from the base case:

\begin{proposition}  Suppose $\Gamma$ is a $d-1$ flag complex where $d$ is $1$ or $2$ and either $\Gamma$ is bipartite or  $\Gamma - e$ is bipartite for some edge $e$ of $\Gamma_i$. Let $V$ be the vertex set of $\Gamma$, $\X = \{ \X_v : v\in V \}$, and let $n = |V|$. Then there is an order $\prec$ on $\X$ and matrices $B$ and $C$ of dimensions $(n-d) \times d$ and $d \times d$, respectively, such that if we let
$$g = \left[ \begin{array}{cc} I_{n-d} & B\\
                                              0 & C\end{array}\right],$$
  then $(\prec, g)$ is a balancing pair for $\Gamma$.    
  
  \begin{proof}  If $d = 0$, take any arbitrary total order $\prec$ on $\X$ and let
  $$g^{-1} = \left[\begin{array}{cc} I_{n-1} &\begin{array}{c}-1\\.\\.\\-1\end{array}\\
                                                            0& 1\end{array}\right],$$
  so that
    $$g = \left[\begin{array}{cc} I_{n-1} &\begin{array}{c}1\\.\\.\\1\end{array}\\
                                                            0& 1\end{array}\right].$$
  The Kind-Kleinschmidt condition is immediate. Our partition of $\X - T_{\prec}$ must simply be $\X_1 = \X - \{x_w\}$, where $x_w$ is the last element of $\X$ with respect to $\prec$. As $\Gamma$ is $0$-dimensional, for any distinct vertices $v_1, v_2$ of $\Gamma$, $x_{v_1}x_{v_2}$ is in $I_{\Gamma}$. If neither $v_1$ or $v_2$ is $w$, then $g(x_{v_1}x_{v_2}) = x_{v_1}x_{v_2}$, so $x_{v_1}x_{v_2} \in \In(gI_{\Gamma})$.
  
  Furthermore, if $v \neq w$, $x_vx_w \in I_{\Gamma}$, and 
  \begin{eqnarray*}
  g(x_vx_w) &= &x_v\sum_{z \in V} x_z\\
                      &=& (\sum_{x_z \prec x_v}x_zx_v) + x_v^2 +( \sum_{x_z \succ x_v}x_vx_z).
                      \end{eqnarray*}
   Then, as $x_zx_v \in gI_{\Gamma}$, for each $x_z \prec x_v$, we have $x_v^2 +( \sum_{x_z \succ x_v}x_vx_z) \in gI_{\Gamma}$, so $x_v^2 \in \In(gI_{\Gamma})$. Thus no monomial in $B_g$ has degree greater than $1$.
   
   If $d =2$ and $\Gamma$ is bipartite, i.e., $2$-colorable, let $V_1$ and $V_2$ be the color classes for $\Gamma$ (for some proper 2-coloring of $\Gamma$). Then identifiying $V_1$ and $V_2$ with the $0$-dimensional complexes on them, $\Gamma$ is a full dimensional subcomplex of $V_1 * V_2$, so the conclusion follows from the $d=1$ case and our lemmas.
  Finally, suppose $d=2$ and $\Gamma$ is not 2-colorable, but $\Gamma-e$ is 2-colorable for some edge $e$ of $\Gamma$. Let $e = \{y, z\}$, and let $A$ and $B$ be the color classes of $\Gamma-e$. Note that $y$ and $z$ must be in the same color class; we may assume that both are in $A$.  Take $\prec$ to be a total order on $\X$ such that the elements of $B$ come before all the elements of $A$, and $y$ and $z$ are the second to last and last elements, respectively, so that $T_{\prec} = \{ y,z\}$. Our partition of $\X -T$ will be $\X_1 = \{x_v : v\in A, v \neq x,y \}$ and $\X_2 = \{x_v : v \in B \}$.
  
  Now, let $C$ be the $(n-2) \times 2$ matrix whose first column is all ones and whose second column has a one in its first $|B|$ rows and zeroes elsewhere. Let
  $$Z= \left[\begin{array}{cc} z_1 & z_2\\
                                                            z_3& z_4\end{array}\right],$$
  and define
  $$g = \left[\begin{array}{cc} I_{n-2} & C\\
                                                            0& Z^{-1}\end{array}\right],$$
so that
  $$g^{-1} = \left[\begin{array}{cc} I_{n-2} & -CZ\\
                                                            0& Z\end{array}\right].$$

We first check the Kind-Kleinschmidt condition. Notice that the rows of $-CZ$ corresponding to variables in $B$ are exactly $\left[\begin{array}{cc} -(z_1 +z_3) & -(z_2 +z_4) \end{array}\right]$, and rows corresponding to variables in $A$ are $\left[\begin{array}{cc} -z_1 & -z_2  \end{array}\right]$.

A facet of $\Gamma$ is  a pair $\{v_1, v_2\}$ where either $v_1 = y$ and $v_2 = z$ or $v_1 \in A$ and $v_2 \in B$. In the first case, the submatrix $g^{-1}$ defined by the intersection of the last two rows with the rows indexed by $v_1$ and $v_2$ is $Z$, in the second it will always have first row $\left[\begin{array}{cc} -(z_1 +z_3) & -(z_2 +z_4) \end{array}\right]$, while the second row will be either $\left[\begin{array}{cc} -z_1 & -z_2  \end{array}\right]$, $\left[\begin{array}{cc} z_1 & z_2  \end{array}\right]$, or  $\left[\begin{array}{cc} z_3 & z_4  \end{array}\right]$. In any case, the rows are linearly independent, so the submatrix has rank 2, and Kind-Kleinschmidt is satisfied.

To complete the proof that $(\prec, g)$ is a balancing pair, it suffices to show that any degree two monomial in $\X_1$ or $\X_2$ lies in $\In(gI_{\Gamma})$.

First, if $x_i \neq x_j$ are both elements of $\X_l$ for $l=1,2$, then they correspond to vertices of the same color class, so $x_ix_j \in I_{\Gamma}$, and $g(x_ix_j)=x_ix_j$, and so $x_ix_j \in \In(gI_{\Gamma})$.

Now suppose $x_i \in X_2$, so $x_i = x_v$ where $v$ is a vertex in $B$. Now, as $\Gamma$ is flag and $1$-dimensional it contains no triangles, so as $\{y, z\} \in \Gamma$, at least one of $\{v,y\}$ or $\{v,z\}$ is not in $\Gamma$. In particular, $x_ix_w \in I_{\Delta}$, where $w$ is either $y$ or $z$. In either case,
$$
g(x_ix_w) = \sum_{j<i}x_jx_i + x_i^2 + S,$$ where $S$ is a sum of degree two monomials occuring later than $x_i^2$ in the revlex order. As $x_jx_i \in gI_{\Gamma}$ for all $j<i$, we then have $x_i^2 + S \in gI_{\Delta}$, so $x_i^2 \in \In(gI_{\Gamma})$.

Finally, suppose $x_i \in \X_1$, so $x_i = x_v$ where $v \in A - \{y,z\}$. Then $x_ix_y$ and $x_ix_z$ are both in $I_{\Delta}$, and so $gI_{\Gamma}$ contains $g(x_ix_y - x_ix_z)$. But
\begin{eqnarray*}
g(x_ix_y - x_ix_z) & = & g(x_ix_y)- g(x_ix_z)\\
			    &=& \sum_{w\in B}x_wx_i + (\sum_{x_j \in \X_2, j<i} x_jx_i) + x_i^2 +S - \sum_{w\in B}x_wx_i + S_2\\
			    & = & x_i^2 + S_1 +S_2,
\end{eqnarray*}
 where $S_1$ and $S_2$ consist of monomials occuring after $x_i^2$ in the revlex order. In particular, $x_i^2 \in \In(gI_{\Gamma})$.
  
\end{proof}  
\end{proposition}

\section{Independence Complexes of Graphs with Large Girth} \label{examples}

Recall that the \emph{independence complex} of a graph $G$ on vertex set $V$ is the simplicial complex $I(G)$ whose faces are exactly the independent sets of $G$, that is, subsets $\tau$ of $V$ such that no two elements of $\tau$ are adjacent in $G$.  A simplicial complex is flag if and only if it is the independence complex of some graph. The aim of this section is to show that Conjecture \ref{kalaiconj}  holds for CM flag complexes arising as independence complexes of graphs of sufficient girth.

Suppose $\Delta = I(G)$ for some graph $G$. Define $\beta(G)$ to be the maximum size of an independent set of $G$, so that $\dim(I(G)) = \beta(G)-1$. If $\Delta$ is Cohen-Macaulay, then $\Delta$ is in particular pure, so all of the maximal independence sets of $G$ have size $\beta(G)$. Such a graph is called \emph{well-covered}. Finbow and Hartnell (see \cite{MR1198396, MR737090}),  gave a characterization of well-covered graphs of large girth:

Let $G$ be a graph on vertex set $V$. A \emph{pendant edge} of $G$ is an edge which is incident to a vertex of degree 1. A \emph{perfect matching} in $G$ is a set of edges $M$ of $G$ such that each vertex of $G$ is in exactly one edge in $M$. 

\begin{theorem} \cite{MR737090} Suppose $G$ is a graph with girth at least $8$. Then $G$ is well-covered if and only if its pendant edges form a perfect matching.
\end{theorem}

It is easy to see that if $G$ is such a graph, then $\beta(G)$ is equal to the number of pendant edges of $G$, or half the number of vertices of $G$.  Then the partition of the vertices of $I(G)$ such that the endpoints of each pendant edge of $G$ constitute a single color class gives a proper $\beta(G)$-coloring of $I(G)$, so $I(G)$ is balanced. Thus Conjecture \ref{kalaiconj} trivially holds for CM-complexes arising in this way.

If we allow smaller girths, however, things become more interesting. Following \cite{MR1198396}, define a 5-cycle in $G$ to be \emph{basic} if it contains no adjacent vertices of degree greater than or equal to 3 in $G$. Let $\PG$ be the set of graphs $G$ such that vertex set of $G$ may be partitioned into two disjoint subsets $P$ and $C$ such that:
\begin{itemize}
\item $P$ contains the vertices in $G$ adjacent to pendant edges of $G$, and the pendant edges form a perfect matching of $P$, and
\item $C$ contains the vertices of the basic 5-cycles of $G$ and the vertices of these 5-cycles give a partition of $C$.
\end{itemize} 

A simple example of a graph in $\PG$ is given in Figure \ref{ex1}.

\begin{figure}[h]
  \begin{center}
    \begin{graph}(6,2)(0,0)
      \graphnodesize{0.1}

     \roundnode{A}(0,1)
      \roundnode{B}(1,2)
      \roundnode{C}(2,1)
      \roundnode{D}(.5,0)
      \roundnode{E}(1.5,0)
       \roundnode{F}(3,1)
      \roundnode{G}(4,2)
      \roundnode{H}(5,1)
      \roundnode{I}(3.5,0)
      \roundnode{J}(4.5,0)
      \roundnode{K}(6,1)
      \roundnode{L}(6,2)

       \edge{A}{B} \edge{B}{C} \edge{A}{D} \edge{D}{E} \edge{E}{C} 
       \edge{F}{G} \edge{G}{H} \edge{F}{I} \edge{I}{J} \edge{J}{H}
       \edge{E}{I}\edge{B}{G}     \edge{H}{K}\edge{K}{L}
        \end{graph}
       
  \end{center}

  \caption{}
  \label{ex1}
\end{figure}

\begin{theorem} \cite{MR1198396} \label{girth} Suppose $G$ is a connected graph of girth greater than or equal to 5. Then $G$ is well-covered if and only if either $G \in \PG$ or $G$ is one of $K_1$, $C_7$, $P_{10}$, $P_{13}$, $P_{14}$, or $Q_{14}$ (see Figure \ref{exceptional}). In the case that $G \in \PG$,  $\beta(G)$ is equal to the number of pendant edges plus twice the number of basic 5-cycles.

\end{theorem}

\begin{figure}[h]
\centering
 \subfigure[$P_{14}$]{
      \begin{graph}(5,4)(.4,0)
         \graphnodesize{0.1}
         
         \roundnode{A}(0,1)
         \roundnode{B}(.7,2.1)
         \roundnode{C}(2,3)
         \roundnode{D}(3.3,2.1)
         \roundnode{E}(4,1)
         \roundnode{F}(2.8,-.1)
         \roundnode{G}(1.2,-.1)
         \roundnode{H}(.9,1.1)
         \roundnode{I}(1.2,1.8)
         \roundnode{J}(2,2.2)
         \roundnode{K}(2.8,1.8)
         \roundnode{L}(3.1,1.1)
         \roundnode{M}(2.4,.5)
         \roundnode{N}(1.6,.5)
         
         \edge{A}{B}\edge{C}{B}\edge{C}{D}\edge{E}{D}\edge{E}{F}\edge{F}{G}\edge{A}{G}\edge{A}{H}\edge{B}{I}\edge{C}{J}\edge{D}{K}\edge{E}{L}\edge{F}{M}
         \edge{G}{N}\edge{I}{K}\edge{K}{M}\edge{M}{H}\edge{H}{J}\edge{J}{L}\edge{L}{N}\edge{N}{I}
      
      \end{graph}

}
\subfigure[$P_{10}$]{
  \begin{graph}(3,5)
         \graphnodesize{0.1}
         
         \roundnode{A}(0,2) 
         \roundnode{B}(1,3)
         \roundnode{C}(2,2)
         \roundnode{D}(2,1)
         \roundnode{E}(0,1)
         \roundnode{F}(4,2)
         \roundnode{G}(4,1)
         \roundnode{H}(3,0)
         \roundnode{I}(5,3)
         \roundnode{J}(5,0)
      \autonodetext{A}[w]{1}  \autonodetext{B}[w]{2}  \autonodetext{C}[n]{3}  \autonodetext{D}[e]{4}  \autonodetext{E}[w]{5}  \autonodetext{F}[e]{6}  \autonodetext{G}[e]{7}  \autonodetext{H}[w]{8}  \autonodetext{I}[e]{9}   \autonodetext{J}[e]{10}
       
         \edge{A}{B}\edge{C}{B}\edge{C}{D}\edge{D}{E}\edge{E}{A}\edge{F}{C}\edge{F}{G}\edge{G}{H}\edge{D}{H}\edge{I}{B}\edge{I}{J}\edge{H}{J}
        
      \end{graph}
}
\subfigure[$P_{13}$]{
  \begin{graph}(6,4)
         \graphnodesize{0.1}
         
         \roundnode{A}(0,2)
         \roundnode{B}(1,3)
         \roundnode{C}(1,1)
         \roundnode{D}(2,3)
         \roundnode{E}(2,2)
         \roundnode{F}(2,1)
         \roundnode{G}(3,4)
         \roundnode{H}(3,0)
         \roundnode{I}(4,3)
         \roundnode{J}(4,2)
         \roundnode{K}(4,1)
         \roundnode{L}(5,3)
         \roundnode{M}(5,1)
         \autonodetext{A}[w]{1}  \autonodetext{B}[nw]{2}  \autonodetext{C}[sw]{3}  \autonodetext{D}[e]{4}  \autonodetext{E}[w]{5}  \autonodetext{F}[e]{6}  \autonodetext{G}[e]{7}  \autonodetext{H}[e]{8}  \autonodetext{I}[n]{9}   \autonodetext{J}[e]{10}\autonodetext{K}[s]{11}\autonodetext{L}[e]{12}\autonodetext{M}[e]{13}

         \bow{A}{G}{.11} \bow{A}{H}{-.11} \edge{B}{D} \edge{C}{B} \edge{D}{E} \edge{E}{F} \edge{F}{C} \edge{D}{G} \edge{F}{H} \edge{E}{J} \edge{G}{I} \edge{I}{J}
          \edge{J}{K} \edge{K}{H} \edge{K}{M} \edge{M}{L} \edge{L}{I}
   \end{graph}       
}
\subfigure[$Q_{13}$]{
  \begin{graph}(5,3)(-.6,0)
         \graphnodesize{0.1}
         
         \roundnode{A}(0,0)
         \roundnode{B}(4,0)
         \roundnode{C}(1,1)
         \roundnode{D}(3,1)
         \roundnode{E}(1,1.5)
         \roundnode{F}(2,1.5)
         \roundnode{G}(3,1.5)
         \roundnode{H}(2,2)
         \roundnode{I}(1,2.5)
         \roundnode{J}(2,2.5)
         \roundnode{K}(3,2.5)
         \roundnode{L}(1,3)
         \roundnode{M}(3,3)
       
         \edge{A}{B}\edge{A}{I}\edge{K}{B}\edge{A}{D}\edge{C}{B}\edge{C}{E}\edge{D}{G}\edge{E}{F}\edge{F}{G}\edge{E}{I}\edge{G}{K}\edge{F}{H}\edge{H}{J}
         \edge{I}{J}\edge{J}{K}\edge{I}{L}\edge{K}{M}\edge{L}{M}
   \end{graph}       
}
 \caption{}
 \label{exceptional}
\end{figure}

Our aim is to show that CM-complexes arising from graphs of girth at least 5 satisfy Conjecture \ref{kalaiconj}. We first address the exceptional cases:

\begin{proposition}\label{notcm} If $G$ is one of $C_7$, $P_{10}$, $P_{13}$, $P_{14}$, or $Q_{14}$, then $I(G)$ is not Cohen-Macaulay.

\begin{proof} It is shown in \cite{MR2302553} that the $I(C_n)$ is Cohen-Macaulay if and only if $n$ is 3 or 5; in particular $I(C_7)$ is not Cohen-Macaulay (this may also be seen by explicitly computing its homology).

Next, note that (using the labeling in Figure \ref{exceptional}) $\lk_{I(P_{10})}(5) = I(C_7)$, so as links in CM complexes are always CM, $I(P_{10})$ is not CM. Similarly, $\lk_{I(P_{13})}(\{10,12\}) = I(C_7)$, so $I(P_{13})$ is not CM.

Finally, we note that $I(P_{14})$ and $I(Q_{13})$ both have dimension 4, but it may be computed (we used the Sage computer algebra system) that each has non-vanishing homology in degree 3. As CM complexes may only have non-vanishing homology in their top degree, neither of these complexes is Cohen-Macaulay. 

\end{proof}

\end{proposition}

We are now in a position to prove the main result of this section:

\begin{theorem} \label{second}Suppose $G$ is a graph of girth at least 5 such that $I(G)$ is Cohen-Macaulay of dimension $d-1$. Then there is a balanced CM complex $\Delta$ of dimension $d-1$ such that $f(I(G)) = f(\Delta)$.

\begin{proof} Suppose that $G$ has connected components $G_1, G_2, \ldots, G_r$. It may easily be seen that $I(G) = I(G_1) * I(G_2) * \cdots * I(G_r)$. It is known  \cite{MR858053} that the join of two complexes is CM if and only if both complexes themselves are CM, in particular each $I(G_i)$ must be CM; furthermore each $G_r$ must have girth at least 5, so by Theorem \ref{girth} and Proposition \ref{notcm} each $G_i$ is either $K_1$ or in $\PG$.

Now, suppose $G_i \in \PG$. Let $\gamma_1, \gamma_2, \ldots, \gamma_j$ be the basic 5-cycles of $G_i$ and $e_1, e_2, \ldots e_l$ the pendant edges of $G$ (so each vertex of $G_i$ is in exactly one $G_s$ or $e_s$). Then $I(G_i)$ is a subcomplex of $I(\gamma_1)*I(\gamma_2) * \cdots * I(\gamma_j) * I(e_1) * \cdots I(e_l)$. Note that each $I(\gamma_s)$ is a 1-complex isomorphic to $C_5$, while each $I(e_s)$ is 0-dimensional. Furthermore, the dimension of $I(G_i) = 2j + l$, so $I(G_i)$ is a full-dimensional subcomplex of this join.

Thus, we see that $I(G)$ is a full-dimensional subcomplex of $\Gamma_1 * \Gamma_2 * \cdots \Gamma_k$ where each $\Gamma_j$ is either 0-dimensional or $C_5$, so we may apply Theorem \ref{main} to complete the proof.

\end{proof}

\end{theorem}

Finally, we note that this class of flag complexes contains a large number of examples:

\begin{corollary} Suppose $G$ is a well-covered graph such that any induced cycle in $G$ has length $5$ and $\beta(G) = d$. Then there is a balanced $(d-1)$-complex $\Delta$ such that $f(I(G)) = f(\Delta)$.

\begin{proof} Clearly the girth of $G$ is either $5$ or $\infty$. In \cite{MR2515394}, Woodroofe showed that if $G$ is well-covered and contains no induced cycles of length other than 5 or 3, then $I(G)$ is CM. Hence we may apply Theorem \ref{second}.

\end{proof}

\end{corollary}

\begin{example} We conclude with an example of a flag complex that does not satisfy the conditions of Theorem \ref{main}. Let $\Delta$ be the flag complex whose $1$-skeleton is the graph $G$ pictured in Figure \ref{badexample}. 

\begin{figure}[h]
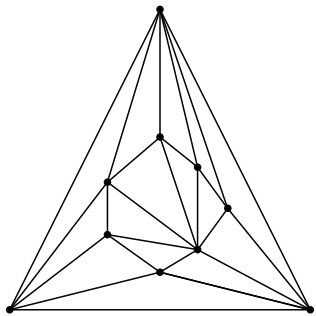

  \begin{center}
    \begin{graph}(4,4)(0,0)
      \graphnodesize{0.1}

     \roundnode{0}(0,0)
     \roundnode{1}(2,4)
     \roundnode{2}(4,0)
     \roundnode{3}(1.3,1)
     \roundnode{4}(1.3,1.7)
     \roundnode{5}(2,2.3)
     \roundnode{6}(2.5,1.9)
     \roundnode{7}(2.9,1.35)
     \roundnode{8}(2.5,.8)
     \roundnode{9}(2,.5)

      \edge{0}{1}\edge{1}{2}\edge{0}{2}\edge{0}{9}\edge{2}{9}\edge{0}{3}\edge{3}{9}\edge{3}{4}\edge{0}{4}\edge{1}{4}\edge{4}{5}\edge{1}{5}\edge{5}{6}\edge{1}{6}\edge{6}{7}\edge{7}{8}\edge{8}{9}\edge{2}{7}\edge{2}{8}\edge{2}{9}\edge{1}{7}\edge{8}{6}\edge{8}{5}\edge{8}{4}\edge{8}{3}
       
        \end{graph}
       
  \end{center}

  \caption{$1$-Skeleton of $\Delta$, $G$}
  \label{badexample}
\end{figure}

Note that $\Delta$ is shellable and hence Cohen-Macaulay (in fact, $\Delta$ is a sphere). The dimension of $\Delta$ is 2.

Suppose $\Delta$ is a full-dimensional subcomplex of some $\Gamma$ of the type described in Theorem \ref{main}. Then $\Gamma$ is of dimension 2, so either $\Gamma  = \Gamma_1 * \Gamma_2 * \Gamma_3$ where each $\Gamma_i$ is 0-dimensional, or $\Gamma = \Gamma_1 * \Gamma_2$ where $\Gamma_1$ is 0-dimensional and $\Gamma_2$ is 1-dimensional. In the former case, it would follow that $G$ is 3-colorable, which it is not.

So suppose  $\Gamma = \Gamma_1 * \Gamma_2$ where $\Gamma_1$ is 0-dimensional and $\Gamma_2$ is 1-dimensional. As $\Delta$ is a pure 2-dimensional subcomplex of $\Gamma$, $\Gamma_1$ must consist of a set of vertices which are pairwise disjoint in $\Delta$ such that every facet of $\Delta$ contains an element of $\Gamma_1$. In other words, $\Gamma_1$ must be an independent set of $G$ that intersects every facet of $\Delta$. One may check by hand that no such independent set exists, and thus $\Delta$ does not satisfy the condition of Theorem \ref{main}.

On the other hand, notice that $f(\Delta) = (1,10,24,16)$, so $h(\Delta) = (1,7,7,1)$. But the simplicial complex $\Omega$ pictured in Figure \ref{ex2} is $3$-colorable and has $f$-vector $(1,7,7,1)$, and therefore, by Theorem \ref{ch}, $(1,7,7,1)$ is the $h$-vector of a balanced Cohen-Macaulay complex of dimension 2. Hence Kalai's conjecture holds for this complex.

\begin{figure}[h] 
 \includegraphics[scale=1]{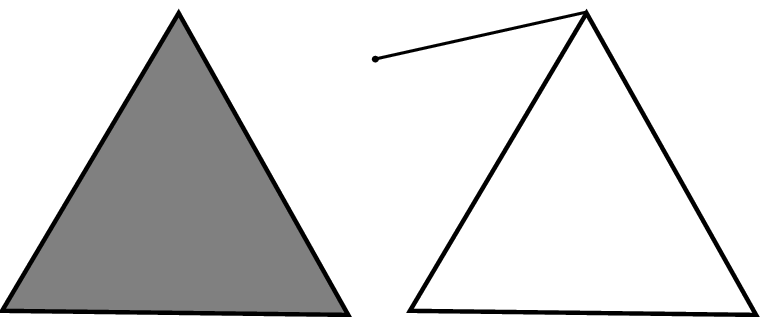}
 \caption{$\Omega$}
 \label{ex2}
\end{figure}

\end{example}

\section{Acknowledgements}
I would like to thank Isabella Novik for many enlightening discussions. This research was partially supported by VIGRE NSF grant DMS-0354131.

\bibliography{jb}

\begin{thebibliography}{10}

\bibitem{MR905148}
A.~Bj{\"o}rner, P.~Frankl, and R.~Stanley.
\newblock The number of faces of balanced {C}ohen-{M}acaulay complexes and a
  generalized {M}acaulay theorem.
\newblock {\em Combinatorica}, 7(1):23--34, 1987.

\bibitem{genbal}
J.~Browder and I.~Novik.
\newblock Face numbers of generalized balanced {C}ohen-{M}acaulay complexes.
\newblock 2009.
\newblock Submitted, preprint available at http://arxiv.org/abs/0909.1134.

\bibitem{MR1322960}
D.~Eisenbud.
\newblock {\em Commutative algebra: With a view toward algebraic geometry},
  volume 150 of {\em Graduate Texts in Mathematics}.
\newblock Springer-Verlag, New York, 1995.

\bibitem{MR1198396}
A.~Finbow, B.~Hartnell, and R.~J. Nowakowski.
\newblock A characterization of well covered graphs of girth {$5$} or greater.
\newblock {\em J. Combin. Theory Ser. B}, 57(1):44--68, 1993.

\bibitem{MR737090}
A.~S. Finbow and B.~L. Hartnell.
\newblock A game related to covering by stars.
\newblock {\em Ars Combin.}, 16(A):189--198, 1983.

\bibitem{MR2302553}
C.~A. Francisco and A.~Van~Tuyl.
\newblock Sequentially {C}ohen-{M}acaulay edge ideals.
\newblock {\em Proc. Amer. Math. Soc.}, 135(8):2327--2337 (electronic), 2007.

\bibitem{MR1018807}
P.~Frankl, Z.~F{\"u}redi, and G.~Kalai.
\newblock Shadows of colored complexes.
\newblock {\em Math. Scand.}, 63(2):169--178, 1988.

\bibitem{MR2391144}
A.~Frohmader.
\newblock Face vectors of flag complexes.
\newblock {\em Israel J. Math.}, 164:153--164, 2008.

\bibitem{MR534824}
B.~Kind and P.~Kleinschmidt.
\newblock Sch\"albare {C}ohen-{M}acauley-{K}omplexe und ihre
  {P}arametrisierung.
\newblock {\em Math. Z.}, 167(2):173--179, 1979.

\bibitem{MR2122283}
I.~Novik.
\newblock On face numbers of manifolds with symmetry.
\newblock {\em Adv. Math.}, 192(1):183--208, 2005.

\bibitem{MR0407036}
G.~A. Reisner.
\newblock Cohen-{M}acaulay quotients of polynomial rings.
\newblock {\em Adv. Math.}, 21(1):30--49, 1976.

\bibitem{MR858053}
C.~Sava.
\newblock On the {S}tanley-{R}eisner ring of a join.
\newblock {\em An. \c Stiin\c t. Univ. Al. I. Cuza Ia\c si Sec\c t. I a Mat.},
  31(2):145--148, 1985.

\bibitem{MR563925}
R.~P. Stanley.
\newblock The number of faces of a simplicial convex polytope.
\newblock {\em Adv. Math.}, 35(3):236--238, 1980.

\bibitem{MR1453579}
R.~P. Stanley.
\newblock {\em Combinatorics and commutative algebra}, volume~41 of {\em
  Progress in Mathematics}.
\newblock Birkh\"auser Boston Inc., Boston, MA, second edition, 1996.

\bibitem{MR0018405}
P.~Tur{\'a}n.
\newblock Eine {E}xtremalaufgabe aus der {G}raphentheorie.
\newblock {\em Mat. Fiz. Lapok}, 48:436--452, 1941.

\bibitem{MR2515394}
R.~Woodroofe.
\newblock Vertex decomposable graphs and obstructions to shellability.
\newblock {\em Proc. Amer. Math. Soc.}, 137(10):3235--3246, 2009.

\end{thebibliography}
\end{document}